\documentclass[11pt]{amsart}
\usepackage{amsthm}
\usepackage{amsmath}
\usepackage{amssymb}
\usepackage{amsfonts}
\usepackage{comment}
\usepackage{graphicx}
\usepackage{hyperref}
\usepackage{cleveref}
\usepackage{leftindex}
\usepackage{longtable}
\usepackage{mathtools}
\usepackage{stmaryrd}
\usepackage{todonotes}
\usepackage{ytableau}

\usepackage[utf8]{inputenc}
\usepackage[margin=1in]{geometry}

\newtheorem{theorem}{Theorem}

\newtheorem{proposition}[theorem]{Proposition}
\newtheorem{lem}[theorem]{Lemma}

\theoremstyle{definition}
\newtheorem{defin}[theorem]{Definition}

\theoremstyle{remark}
\newtheorem{rmk}[theorem]{Remark}

\newcommand{\cone}{\mathrm{cone}}
\newcommand{\cx}{\mathrm{cx}}

\newcommand{\Span}{\mathrm{span}}

\newcommand{\Oh}{\mathcal{O}}
\newcommand{\R}{\mathbb{R}}
\newcommand{\T}{\mathcal{T}}

\newcommand{\cA}{\mathcal{A}}

\newcommand{\cM}{\mathcal{M}}

\newcommand{\cS}{\mathcal{S}}

\newcommand{\flip}{\mathrm{flip}}
\newcommand{\RRS}{\mathrm{RRS}}

\begin{document}
\title{Oriented matroid structures on rank 3 root systems}
\author{Grant Barkley}
\address{Department of Mathematics, Harvard University, Cambridge, MA}
\email{gbarkley@math.harvard.edu}

\author{Katherine Tung}
\address{Department of Mathematics, Harvard University, Cambridge, MA}
\email{katherinetung@college.harvard.edu}

\begin{abstract}
    We show that, given a rank 3 affine root system $\Phi$ with Weyl group $W$, there is a unique oriented matroid structure on $\Phi$ which is $W$-equivariant and restricts to the usual matroid structure on rank 2 subsystems. Such oriented matroids were called oriented matroid root systems in \cite{DyerWang}, and are known to be non-unique in higher rank. We also show uniqueness for any finite root system or ``clean'' rank 3 root system (which conjecturally includes all rank 3 root systems).  
\end{abstract}
\maketitle

\section{Introduction} \label{sec:intro}

In this note, we examine the possible oriented matroid structures on a root system.
In general, there may be combinatorially distinct root systems associated to the same Coxeter group $W$ (in the sense that their underlying oriented matroids are distinct) \cite{DyerWang}. In contrast, for rank $3$ Coxeter groups there is a unique oriented matroid underlying any realization of its root system.
Dyer and Wang conjecture \cite[Section 6.8]{DyerWang} that this oriented matroid is the unique oriented matroid structure on the root system with $W$-invariant circuits, so that no set of positive roots forms a circuit and so that it has the same one-, two-, and three-element circuits as any realization. Such a structure is called an \emph{oriented matroid root system}. Our main result is a proof of this conjecture for rank 3 root systems which are \emph{clean} (see Definition \ref{def:clean}). This property was studied in \cite{barkley2024affineextendedweakorder}, where it was conjectured that all rank 3 root systems are clean. 

\begin{theorem}\label{thm:rank3}
    Let $\Phi$ be a clean rank 3 root system. Then $\Phi$ has a unique oriented matroid root system structure.
\end{theorem}

It was proven in \cite{barkley2024affineextendedweakorder} that all rank 3 affine root systems are clean, so our result is unconditional when $\Phi$ is an affine root system. We also show

\begin{theorem}\label{thm:finite}
    Let $\Phi$ be a finite root system. Then $\Phi$ has a unique oriented matroid root system structure.
\end{theorem}

\section{Preliminaries} \label{sec:prelim}

\subsection{Dictionary between hyperplane arrangement and oriented matroid language} \label{subsec:dict}
This table gives some heuristics for translating between hyperplane arrangement language and the corresponding oriented matroid language, which we will use in the proof below.

 \begin{longtable}[c]{| c | c |}
 
 \hline
 $\cA$ is a hyperplane arrangement &  $\cM = (E,\ast,\cx)$ is an oriented matroid \\
 \hline
 Unit normal vector to a hyperplane & Element of $E$ \\
 Half-space bounded by a hyperplane & Element of $E$ \\
 Negation of a half-space or normal vector & $e\mapsto e^*$ \\
 Regions of $\cA$ & Topes of $\cM$ \\
 {Subarrangement} & {$\ast$-closed subset of $E$} \\
 Adjacent regions are separated by 1 hyperplane & Reorientation property satisfied \\
 $\alpha$ is a non-negative combination of vectors in $S$ & $\alpha\in \cx(S)$ \\

Deletion $\cA\setminus \{H\}$ & Restriction $\cM\restriction E\setminus\{e,e^*\}$\\
 \hline
 \endfirsthead
\end{longtable}

\subsection{Oriented matroids}\label{subsec:oriented-matroid}
\begin{defin}\label{def:om-dyer} \cite{DyerWang}
    An \emph{oriented matroid} is a triple $\cM=(E, *, \cx)$ where $E$ is a set with a fixed-point free involution map $*: E \to E$ and $\cx$ a closure operator on $E$ such that 
    \begin{enumerate}
        \item $\cx(X)^\ast = \cx(X^\ast)$,
        \item If $x \in \cx(X \cup \{x^*\})$, then $x \in \cx(X)$,
        \item If $x \in \cx(X \cup \{y^*\})$ and $x \not\in \cx(X)$ then $y \in \cx(X \setminus \{y\} \cup \{x^*\}),$
        \item If $x \in \cx(X)$, then there exists a finite set $Y \subseteq X$ such that $x \in \cx(Y)$.
    \end{enumerate}
\end{defin}

Let $\Gamma$ be a subset of some real, finite-dimensional vector space $V$ which does not contain $0$. Let $\cone(A) := \{\sum_{i \in I} k_i v_i \mid v_i \in A \cup \{0\}, k_i \in \R_{\ge 0}, |I| < \infty\}.$ Let $\cone_{\Gamma}(A) := \cone(A) \cap \Gamma$, where $A \subseteq \Gamma.$
\begin{defin}
    An oriented matroid is called \emph{realizable} if it can be written as $(\Gamma, -, \cone_\Gamma)$ where $\Gamma$ is a nonempty set satisfying $\textbf{0} \not\in \Gamma$ and $\Gamma = -\Gamma.$
\end{defin}

\begin{rmk}
    Oriented matroids are often presented using circuits, which are collections of signed sets. Our definition does not use signed sets. We explain briefly how to convert between these cryptomorphic conventions (see also \cite[Exercises 3.9-3.13]{Bjorner} and the references there). Let $E_0$ be the ground set for an oriented matroid presented using signed sets. Then the ground set using our definition is $E\coloneqq E_0 \sqcup E_0^\ast$, where the elements of $E_0^\ast$ are in bijection with $E_0$ and can be thought of as formal symbols $e^*$ for each $e\in E_0$. We extend this to an involution on $E$ by defining $e^{**} = e$. 

    To define the closure operator $\cx$ on $E$, it is sufficient to describe which sets are closed (i.e. satisfy $\cx(X) = X$). We say a set $X$ is closed with respect to a circuit $\Oh$ if, for all $e\in \Oh^+$,
    \[ \Oh^+\setminus\{e\} \cup \{ e^* \mid e \in \Oh^-\} \subseteq X \implies e^*\in X,    \]
    and similarly for all $e\in \Oh^-$. The sets which are closed under $\cx$ are those which are closed with respect to all circuits, and this determines the oriented matroid $(E,\ast, \cx)$. We will not use signed sets in the remainder of the note.
\end{rmk}

\begin{defin}
An element contained in $\cx(\varnothing)$ is called a \emph{loop}, and an oriented matroid without any loops is called \emph{loopless}. A loopless oriented matroid such that $\cx(\{x\}) = \{x\}$ for all $x\in E$ is said to be \emph{reduced} (or \emph{simple}). 
\end{defin}
Every oriented matroid $\cM = (E,\ast,\cx)$ has an associated reduced oriented matroid, its \emph{reduction}. Define an equivalence relation $\sim$ on $E\setminus \cx(\varnothing)$ so that the $\sim$-equivalence class containing $x$ is $\cx(x)\setminus \cx(\varnothing)$. Define $\overline{E} = (E\setminus \cx(\varnothing))/\sim$. The involution and closure operators on $E$ induce corresponding operators on $\overline{E}$. Then $(\overline{E} , \ast, \cx)$ is an oriented matroid, which we call the \emph{reduction} of $\cM$ and denote by $\overline{\cM}$.

For simplicity, we will always assume our oriented matroids are reduced in the rest of the paper. In particular, when we work with realizable oriented matroids $(\Gamma,-,\cone_\Gamma)$, then we will replace $\Gamma$ with $\overline{\Gamma}$ if necessary. The resulting matroid is still realizable, since we get an equivalent oriented matroid by picking one representative for each ray in $\Gamma$ in a way that guarantees they are closed under involution.

\begin{rmk}\label{rmk:rrs-omrs}
Let $\Gamma=\Phi$ is a root system with a base $\Pi$, associated to a Coxeter system $(W,S)$, and let $T$ be the set of conjugates of elements of $S$. Then there is a canonical bijection $\overline{\Phi} \cong T\times \{\pm 1\}$ \cite{DyerWang}. In particular, any two root systems for $(W,S)$ are in canonical bijection after reduction, allowing us to compare the oriented matroid structures arising from different root systems with the same Coxeter group. Such an oriented matroid is called a \emph{reduced root system} (RRS) oriented matroid.
\end{rmk}

\begin{defin}
A subset $X\subset E$ is \emph{asymmetric} if $X^*\cap X = \varnothing$. A subset which is maximal among $\cx$-closed, asymmetric sets is called a \emph{tope}. 
\end{defin}

\begin{rmk}
    If $\cM=(\Gamma,-,\cone_\Gamma)$ is a finite realized oriented matroid, then the topes of $\cM$ biject with the regions of the hyperplane arrangement dual to $\Gamma$. If $\cM$ is an infinite realized oriented matroid, then the topes biject with the \emph{weakly separable sets} in the hyperplane arrangement (see \cite{barkley2024combinatorialdescriptionsbiclosedsets}).
\end{rmk}

If $(E,\ast,\cx)$ is a reduced oriented matroid, then for any $x\in E$ and any tope $R$, it is true that $|\{x,x^*\}\cap R| = 1$. The subsets $R\subseteq E$ such that $|\{x,x^*\}\cap R| = 1 $ for all $x\in E$ are called \emph{topal}. 
More generally, if $(E,\ast,\cx)$ is a reduced oriented matroid, and $\T$ denotes the collection of its topes, then the tuple $(E,\ast,\T)$ uniquely determines $\cx$ \cite[Exercise 3.30]{Bjorner}. 

Given sets $S_1,S_2$, let $S_1 \Delta S_2$ denote their symmetric difference. Given a topal set $R$ and $x\in E$, the \emph{flip} of $R$ across $x$ is the topal set $\flip_x(R) \coloneqq R \Delta \{x,x^*\}$. 

\begin{defin}
We define the \emph{separating set} between two topal sets to be their symmetric difference, modulo $\ast$:
\[ \cS(R_1,R_2) \coloneqq (R_1 \Delta R_2)/\ast. \]

Two topal sets are said to be \emph{adjacent} if $|\cS(R_1,R_2)|=1$. 
The topal set $R_2$ is said to be \emph{between} topal sets $R_1$ and $R_3$ if $\cS(R_1,R_3) = \cS(R_1,R_2) \sqcup \cS(R_2,R_3)$.
\end{defin}

\begin{proposition}[Reorientation property \cite{Handa}]
        Let $\cM = (E,\ast, \cx)$ be a reduced oriented matroid.
        For all topes $R_1$ and $R_3$, if there are no topes strictly between $R_1$ and $R_3$, then $R_1$ and $R_3$ are adjacent.
\end{proposition}

We now define the restriction of an oriented matroid $\cM$ to a subset of its ground set (sometimes called \emph{deletion}). 

\begin{defin}\label{def:restriction}
    Let $\cM = (E,\ast, \cx)$ be an oriented matroid, and let $E'\subseteq E$ be a subset so that $\ast E' = E'$. Then the \emph{restriction} of $\cM$ to $E'$, denoted $\cM\restriction E'$, is the oriented matroid $(E',\ast, \cx_{E'})$, where $\cx_{E'}$ is the closure operator on subsets of $E'$ defined by
    \[ \cx_{E'}(X) \coloneqq \cx(X)\cap E'. \]
\end{defin}

\subsection{Biclosed sets and oriented matroid root systems}\label{subsec:omrs}

Here we introduce a closure operator which is not usually an oriented matroid closure. It is denoted $c_2$ because it is the ``rank $2$ convex closure.'' This operator is the same as the one used in the definition of closed, coclosed, and biclosed in \cite{barkley2024affineextendedweakorder}. 
\begin{defin} \cite{barkley2024affineextendedweakorder}
Let $X$ be a subset of some finite dimensional vector space $V$, and let $B \subseteq X.$ Then
    \begin{enumerate}
        \item $B$ is $c_2$-\emph{closed} in $X$ if, whenever $\alpha,\beta \in B$ and $\gamma \in \Span_+(\alpha, \beta) \cap X$, then $\gamma \in B.$ 
        \item $B$ is $c_2$-\emph{coclosed} in $X$ if $X \setminus B$ is closed in $X.$
        \item $B$ is $c_2$-\emph{biclosed} in $X$ if $B$ is closed and coclosed in $X.$
    \end{enumerate}
    Let $R$ be a topal subset of an oriented matroid of the form $(\overline{\Phi}, \ast, \cx)$ for a root system $\Phi$. Then
    \begin{enumerate}
        \setcounter{enumi}{3}
        \item $R$ is a \emph{quasitope} if $\{\alpha\in \Phi^+ \mid \overline{\alpha} \in X \}$ is $c_2$-biclosed in $\Phi^+$.
    \end{enumerate}
\end{defin}

We note that the map sending a quasitope $R$ to the $c_2$-biclosed set $\{\alpha\in \Phi^+\mid \overline{\alpha}\in X\}$ is a bijection between quasitopes of $\overline{\Phi}$ and $c_2$-biclosed sets of $\Phi^+$. The inverse map sends a $c_2$-biclosed set $B$ to $\overline{B} \cup -(\overline{\Phi}^+\setminus \overline{B}).$

Fix a Coxeter system $(W,S)$. Let $\Phi$ be a root system for $W$ and consider its reduction $\overline{\Phi}$. (Alternatively, we may identify $\overline{\Phi}$ with $T\times \{\pm 1\}$ \textit{vis-\`a-vis} Remark \ref{rmk:rrs-omrs}). 
    
\begin{defin} \label{def:omrs}
An \emph{oriented matroid root system} (OMRS) for $W$ is a reduced oriented matroid $(\overline{\Phi}, \ast, d)$ satisfying the following properties:
\begin{itemize}
    \item $d$ is $W$-equivariant: if $X\subseteq \overline\Phi$ and $w\in W$, then $d(wX) = w\, d(X)$,
    \item $\overline\Phi^+$ is $d$-closed, and
    \item If $X \subseteq \overline{\Phi}$ has span with dimension at most $2$, then $d(X) = \cone_{\Phi}(X)$. 
\end{itemize}
\end{defin}
    The third condition in Definition \ref{def:omrs} has two hidden assumptions. Since a Coxeter group may have multiple different root system realizations, the definition implicitly assumes that the choice of realization for the root system does not affect $\cone_\Phi(X)$ when $X$ has span of dimension at most 2. Furthermore, it assumes that the condition that $X$ has span of dimension at most 2 itself does not depend on the realization of $\Phi$. These assumptions are justified in \cite{DyerWang}. 

The following lemma is straightforward.
\begin{lem}\label{lem:topeisquasitope}
    Every tope of an OMRS is a quasitope.
\end{lem}

We have seen each root system associated to $W$ gives rise to a reduced root system (RRS) oriented matroid $(\overline{\Phi},\ast,\cone_\Phi)$. Each RRS oriented matroid is an OMRS, so in particular each of its topes is a quasitope. A root system is clean if the converse also holds.

\begin{defin}[\cite{barkley2024affineextendedweakorder}]
\label{def:clean}
    Let $\cM= (\overline{\Phi}, \ast, d)$ be the RRS OMRS for the root system $\Phi$. We say that $\Phi$ is \emph{clean} if every quasitope of $\cM$ is a tope.
\end{defin}

\section{Proofs}\label{sec:results}

\begin{proof}[Proof of \Cref{thm:rank3}]
    Let $\cM=(\overline{\Phi},\ast,d)$ be an OMRS on a clean rank $3$ root system $\Phi$. We will show that $\cM$ coincides with $\cM_{\RRS}$, the RRS oriented matroid for $\Phi$.
    
    Let $\beta_1,\beta_2,\ldots \in \Phi^+$ be some ordering of the roots.  We will show by induction on $n$ that for the order ideal $I\coloneqq \{\beta_1,\ldots, \beta_n\}$, the restriction $\cM\restriction \pm I$ coincides with $\cM_{\RRS}\restriction \pm I$. Since oriented matroids are determined by their topes, it is enough to show that the two oriented matroids have the same topes. We observe that every tope of $\cM\restriction \pm I$ is a tope of $\cM_{\RRS}\restriction \pm I$, because every tope of $\cM\restriction \pm I$ is the restriction of a tope of $\cM$, every tope of an oriented matroid is a quasitope, any two OMRS structures on $\Phi$ have the same quasitopes, and the quasitopes of $\cM_{\RRS}$ coincide with the topes (by cleanliness). 
    Hence, it is enough to show that every tope of $\cM_{\RRS}\restriction \pm I$ is a tope of $\cM \restriction \pm I$.  
    
    Topes in $\cM_{\RRS}\restriction \pm I$ can be identified with regions of the hyperplane arrangement dual to $\{\beta_1,\ldots, \beta_n\}$. From here on, we refer to such topes as ``regions'', and reserve the word ``tope'' for topes of $\cM\restriction \pm I$. In this language, our goal is to show that every region is a tope. 
    
    We can represent the order ideal $I$ pictorially by the following system. 
   We draw the hyperplane arrangement dual to $\{\beta_1,\dots,\beta_n\}$ in $\R^3$ and then intersect with a plane to get a 2D cross-section like the one in Figure \ref{fig:figure-b}. We require this plane to intersect each of the hyperplanes $\beta_1^\perp, \dots, \beta_n^\perp$ in a distinct line. We observe:
    \begin{enumerate}
       \item Regions in a cross-section are a subset of the regions of the arrangement.

        \item A region may appear differently in different cross-sections: for example, in one cross-section it may be bounded on all sides whereas in another it may be unbounded.
        \item Any region in an essential central hyperplane arrangement is contained in some half-space, and thus there exists a cross-section where it appears as (convex and) bounded on all sides.

        \item If a region is bounded in a cross-section, then its adjacent regions also also appear as regions of that cross-section
    \end{enumerate}

    For convenience, from now on we refer to the cross-sections of hyperplanes $\beta_i^\perp$ also as $\beta_i^\perp$. 

    If $I=\{\beta_1,\ldots, \beta_n\}$ spans a subspace of dimension at most 2, then $\cM\restriction \pm I$ and $\cM_{\RRS}\restriction \pm I$ coincide by definition of an OMRS. So assume that $\{\beta_1,\ldots, \beta_n\}$ spans a subspace of dimension at least 3.

    Our inductive step is mostly reliant on the following implications of the inductive hypothesis, which asserts that all regions of $\{\beta_1^\perp,\ldots,\beta_{n-1}^\perp\}$ are topes.
    \begin{enumerate}
        \item In a cross-section, any region not bordering $\beta_n^\perp$ is a tope.
        \item Suppose that $\beta_n^\perp$ bisects some region of $\{\beta_1^\perp,\ldots, \beta_{n-1}^\perp\}.$ 
        Then the two new regions created have the property that at least one is a tope.
    \end{enumerate}
    We also rely on the following fact, equivalent to the reorientation property for the topes of a finite oriented matroid. Call a sequence of topes a \textit{gallery} if consecutive topes are adjacent. Then there is a gallery between any two topes $R_1,R_2$ so that each hyperplane in $\cS(R_1,R_2)$ is crossed exactly once. Call such a gallery \textit{minimal.}

    Now, assume that $R$ is a region of $\{\beta_1^\perp, \ldots, \beta_n^\perp\}$ which is not a tope. By the assumption on rank, the region $R$ has at least $3$ walls.
   
    We claim that if $R$ has $3$ walls, then the oriented matroid $\cM$ will view the walls as linearly dependent. Indeed, say the three walls are $\alpha^\perp,\beta^\perp,\gamma^\perp$. Consider the restriction $\cM\restriction \pm\{\alpha,\beta,\gamma\}$. Then the restriction of $R$ to this subset is not a tope, which means that there are less than $8$ topes of $\cM\restriction\pm\{\alpha,\beta,\gamma\}$. This can occur only if $\gamma \in d(\{\pm\alpha,\pm\beta\})$. However, this contradicts the definition of an OMRS, since $\gamma \not\in \cone_\Phi(\{\pm\alpha,\pm\beta\})$. We conclude that $R$ must have at least $4$ walls.

    If $R$ has $N\geq 4$ walls, then there exists a cross-section where it appears as a closed polygon with $N$ sides. At least one of those sides is contained in $\beta_n^\perp$, by observation (1). Consider Figure \ref{fig:figure-b}, which depicts the hyperplanes bordering $R$. Let the region obtained by starting from $R$ and crossing $\beta_n^\perp$ be $T.$ Label the vertices of $R$ by $v_1, v_2, \dots, v_N$ counterclockwise so that $v_1v_2$ is the side contained in $\beta_n^\perp.$ Let $i$ be chosen so that $v_3v_4$ is the side on $\beta_i^\perp$. 
    Let $X$ be the region obtained by starting from $R$ and crossing $\beta_i^\perp.$ Then $X$ does not have $\beta_n^\perp$ as a wall. Hence $X$ and $T$ are both topes by observations (1) and (2), respectively. Since $\cS(X,T)=\{\pm \beta_i, \pm \beta_n\}$ and $\beta_n^\perp$ is not a wall of $X$,
    any minimal gallery from $X$ to $T$ must pass through $R$, so $R$ is a tope. 
\begin{figure}

\tikzset{every picture/.style={line width=0.75pt}} 

\begin{tikzpicture}[x=0.75pt,y=0.75pt,yscale=-1.25,xscale=1.25]

\draw (100,100) -- (200,200);

\draw(200,100) -- (300,200);

\draw (270,100) -- (61,193.5);

\draw (350,100) -- (141,193.5);

\draw(100,130) -- (400,130);

\node[below right] at (252,142) {$v_1$};
\node[above left] at (238,145) {$v_2$};
\node[above left] at (210,130) {$v_3$};
\node[below left] at (148,142) {$v_4$};
\node[above left] at (186,194) {$v_5$};

\node at (200,150) {$R$};

\node[anchor=west] at (148,140) {$X$};

\node[anchor=west] at (240,138) {$T$};

\node at (200,90) {$\beta_{n}^\perp$};

\node at (275,90) {$\beta_{i}^\perp$};

\end{tikzpicture}
 \caption{Configuration of roots in $I$ for $R$ having least 4 sides.}
    \label{fig:figure-b}
\end{figure}
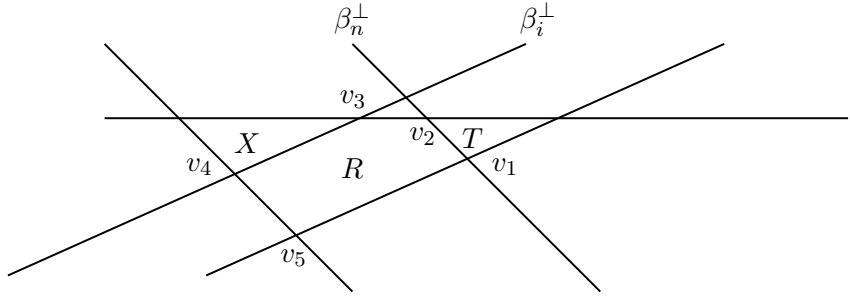
\end{proof}

\begin{rmk}
    Effectively, the proof of Theorem \ref{thm:rank3} shows that one cannot have a rank $3$ oriented matroid whose topes are a proper subset of the topes of a realizable rank $3$ oriented matroid.
\end{rmk}

\begin{proof}[Proof of \Cref{thm:finite}]
     Let $\Phi$ be a finite root system, and let $\cM=(\overline{\Phi}, \ast, d)$ be an OMRS on $\Phi$.
     First, we note that $\overline{\Phi}^+$ is a tope because it is a $d$-closed topal set.
     Using $W$-equivariance of $d$, we can generate a set of $|W|$ distinct topes $\{w\overline{\Phi}^+ \mid w \in W\}$. Each of these topes is a quasitope, and the number of quasitopes in a finite root system is exactly $|W|$ \cite[Lemma 4.1]{Dyerextendedweakorder}. Hence the topes of the OMRS are exactly the quasitopes, and we conclude that the OMRS structure is unique since an oriented matroid is determined by its topes.
    
\end{proof}

\section*{Acknowledgments}
The authors thank Matthew Dyer and David Speyer for helpful conversations.
Author Grant Barkley was supported by NSF grant DMS-1854512 and author Katherine Tung was supported by NSF grant DMS-2053288 during the time this research was conducted.

\bibliographystyle{plain}
\bibliography{main}
\end{document}